\documentclass[12pt]{amsart}
\usepackage{amsthm,amsmath,amssymb}
\usepackage{addfont}
\usepackage{mathrsfs}
\usepackage[T1]{fontenc}
\numberwithin{equation}{section}

\def\cai{{\mathcal I}}
\def\cj{{\mathcal J}}

\def\cp{{\mathcal P}}

\def\cas{{\mathcal S}}

\def\bc{{\mathbb C}}

\def\br{{\mathbb R}}

\def\d{\delta}

\def\g{\gamma}

\def\l{\lambda}

\def\X{\Xi}

\theoremstyle{plain}
\newtheorem{lemma}{Lemma}[section]
\newtheorem{proposition}[lemma]{Proposition}
\newtheorem{theorem}[lemma]{Theorem}
\newtheorem{corollary}[lemma]{Corollary}
\theoremstyle{definition}

\newtheorem{remark}[lemma]{Remark}
\newtheorem{definition}[lemma]{Definition}
\newtheorem{example}[lemma]{Example}

\begin{document}

\title[The Schur multiplier norm and its dual norm]{\textsc{The Schur multiplier norm and  its dual norm}}

\author[E.~Christensen]{Erik Christensen}
\address{\hskip-\parindent
Erik Christensen, Mathematics Institute, University of Copenhagen, Copenhagen, Denmark.}
\email{echris@math.ku.dk}
\date{\today}

\begin{abstract}   For a complex self-adjoint  $n \times n $ matrix $A$ we show that its Schur multiplier norm is determined by $$ \|A\|_S = \min \{\, \|\mathrm{diag}(P)\|_\infty \, :\,   - P \leq A \leq P \, \}.$$
 For a norm $\||.\||$ on $M_{(m,n)}(\bc)$ we define for an $X$ in $M_{(m,n)}(\bc)$  the
  dual norm, $\||X\||^*,$ of $X$ by the formula: 
  $$ \||X\||^*:=
  \sup \{ |\mathrm{Tr}_n(Y^*X)| : Y \in M_{(m,n)}(\bc), \, \||Y\|| \leq 1\}.$$ For a vector $\lambda$  in $\bc^n $ we let $\Delta(\lambda)$ denote the diagonal $n \times n$ matrix with entries from $\lambda. $ 
 The dual norm of the Schur multiplier norm on $M_n(\bc) $ is for a self-adjoint
  matrix $A$ given by the formula:  $$ \|A\|_S^*  = \min \{ \, \mathrm{Tr}_n\big(\Delta(\lambda)\big)\, :\, \lambda \in \br^n,  \, - \Delta(\lambda) \leq  A \leq \Delta(\lambda)\,\}. $$  
  We study this minimization problem as a {\em formal} linear program.  
 \end{abstract}

\maketitle

\section{Introduction}

Issai Schur studied in  \cite{Sc} the entry wise product of two scalar  $ n \times n $ matrices $(X \circ Y)_{(i,j)}:= X_{(i,j)}Y_{(i,j)}$ and he looked at the multiplier norm - now the Schur multiplier norm $\|X\|_S$ -   where the  factor $X$ is fixed and $Y$  in $M_n(\bc)$ is equipped with the operator norm, $\|Y\|_\infty.$   The following theorem is a reformulation  of {\em Satz V} in Schur's article. 
\begin{theorem} \label{Schur}
Let $X$  be a positive complex $n \times n$ matrix, then \newline $\|X\|_S = \|\mathrm{diag}(X)\|_\infty.$ 
\end{theorem} 
This is a remarkably simple formula, and here we  present, in Theorem \ref{SchurNorm}, an extension of this formula to the case of self-adjoint matrices. This result is mentioned in the abstract. The new formula is based on the concept of order on self-adjoint matrices and norms of diagonals, and in the case of a positive matrix we just get Schur's theorem back. 
On the other hand our result is very closely related  to previous results by Haagerup \cite{UH}, which here is named Theorem \ref{Haa}, and an extension of his results by Ando and Okubo \cite{AO}, but we think, that our version of their findings is closer  to the original formulation of Schur's theorem.
At the top of page 112 in Paulsen's book \cite{Pa}, he states a version of the result we present here, but he references a non existent Exercise for a proof. The item (v) of his  Exercise 8.8 offers a proof of Haagerup's result \cite{UH}.  

In the abstract we defined the dual norm of a norm on $M_{(m,n)}(\bc)$ and our definition agrees with the one  R. Mathias uses in his study \cite{Ma1} of a pair of dual norms on $M_{(m,n)}(\bc). $ We want to thank the referee, who directed us to Mathias' work. In our first version of this article we were not able to benefit from Mathias' results, which have simplified our exposition a lot. In the first version of this article we used a result from \cite{C3}, which shows that the dual norm to the Schur multiplier norm is the norm, which we have denoted $\|.\|_{cbB}. $ The reason for this name is that for a complex $ m \times n $ matrix $ X$ the norm $\|X\|_{cbB}$ denotes the completely bounded norm of the bilinear form $B_X$ on $(\bc^m, \|.\|_\infty) \times ( \bc^n, \|.\|_\infty) $ with kernel $X.$ This concept is defined in  \cite{CS} and some of its  properties are established in \cite{C3, C4}. We will need some results from \cite{C3, C4}, but the general theory of completely  bounded bilinear forms, as described in \cite{Pa},  is not needed in this article, so we will base our construction of dual norms on  \cite{Ma1} and use the symbol $\|.\|_{\cas^*} $ for the dual norm to the Schur multiplier norm instead of the symbol $\|.\|_{cbB},$ which was used in the preprint. 

The second formula mentioned in the abstract was proven in  Theorem 3.5 of \cite{C4} for positive matrices and it turned out, while realizing how Schur's result could be generalized to self-adjoint matrices, that a similar extension also was possible in the dual case. The duality just mentioned may be seen as a shadow of the well known duality between an $\ell^\infty-$norm, played by the $\|.\|_S$ and an $\ell^1-$norm played by $\|.\|_{\cas^*}.$

Section \ref{Exa} contains an example which shows that for a complex self-adjoint $2n \times 2n$ matrix $Y$ we may have $$ - \sqrt{Y^2} \leq Y \leq \sqrt{Y^2} \text{ and } \|\mathrm{diag}(\sqrt{Y^2})\|_\infty = \sqrt{n} \|Y\|_S,$$ so the obvious candidate for the optimal positive matrix in Theorem  \ref{SchurNorm}, when applied to this $Y$ does not work very well. 
This example is also related to Walter's inequality from \cite{Wa}.

In Section \ref{LP} we transform the statement in Corollary  \ref{cbBSa} into a question of solving an {\em infinite formal version of a linear program} and then we show that the corresponding infinite formal  dual program has the same optimal value and its optimal solutions are closely connected to the duality between the $\cas^*$-norm and the Schur multiplier norm.  
R. A. Horn uses the words {\em Hadamard multiplier} instead of the words {\em Schur multiplier,} and we recommend his survey article \cite{Ho} as an introduction to the literature on Hadamard multipliers.

Finally we will like to thank both of the referees, who have read the preprint very carefully and been most helpful by  suggesting improvements in the presentation, detected errors in the arguments  and directed us to look into Mathias' works \cite{Ma1, Ma2}.   

\section{ The Schur multiplier norm of a self-adjoint matrix} 

In \cite{UH} Haagerup presents a  characterization, of the Schur multiplier norm of a  complex matrix, which we recall in the following  theorem. The proof of this theorem also appears in Ando and Okubo's article \cite{AO}, and as  item (v) of Exercise 8.8 in Paulsen's book \cite{Pa}. The story behind this result is long and the standard form of a matrix with Schur multiplier norm at most one, can be found in Paulsen's book as Theorem 8.7. 

\begin{theorem} \label{Haa} 
Let $A$ be a complex $m \times n $ matrix, then the Schur multiplier norm of $A$ is at most 1 if and only if there exist positive matrices $P_1 $ in $M_m(\bc)$ and $P_2 $ in $M_n(\bc)$ such that 
\begin{equation}
\begin{pmatrix}
P_1 & A\\ A^* & P_2
\end{pmatrix} \geq 0 \text{ and  for } i \in \{1,2\}:\, \, \|\mathrm{diag}(P_i)\|_\infty \leq 1.
\end{equation}
\end{theorem}

This theorem may take a form which avoids block matrices in the case when $A$ is a self-djoint matrix in $M_n(\bc).$ 

\begin{theorem} \label{SchurNorm}
Let $A$ be a self-adjoint matrix in $M_n(\bc),$  then the Schur multiplier norm of $A$ is at most 1 if and only if there exists a positive matrix $P$ in $M_n(\bc)$ such that 
\begin{equation}
-P \leq A \leq P \text{ and  } \| \mathrm{diag}(P)\|_\infty \leq 1.
\end{equation}
\end{theorem}
\begin{proof}
The arguments  in  this proof and  in the proof of the following corollary are based on Schur's basic theorem, which tells  that for a positive matrix $P$ we have the identity $\|P\|_S = \|\mathrm{diag}(P)\|_\infty.$ 

Let us  first assume that $\|A\|_S \leq 1 $ and let positive matrices $P_1$ and $P_2$ be found in $M_n(\bc)$ such that $$ \begin{pmatrix}
P_1& A\\A & P_2
\end{pmatrix}  \geq 0 \text{ and } \|\mathrm{diag}(P_i)\|_\infty \leq 1.$$
Let $P:= (P_1 + P_2)/2,$ then $P$ is positive and $\| \mathrm{diag}(P) \|_\infty \leq 1.$ Further we define a {\em positive block matrix}  $Q$ by 
$$Q := \begin{pmatrix}
P & A \\A & P
\end{pmatrix} = \frac{1}{2} \begin{pmatrix}
P_1& A \\A & P_2
\end{pmatrix} + \frac{1}{2}\begin{pmatrix}
0 & I_n\\I_n & 0
\end{pmatrix} \begin{pmatrix}
P_1& A \\A & P_2
\end{pmatrix}\begin{pmatrix}
0 & I_n \\I_n & 0
\end{pmatrix}  .$$
 The block matrices $M_2\big(M_n(\bc)\big)$  may be identified with the tensor product $M_n(\bc)\otimes M_2(\bc)$ and in this context we define 2 orthogonal projection matrices $E $ and  $F$ in $M_2(\bc)$ by $$ E = \begin{pmatrix}
1/2 & 1/2\\ 1/2 & 1/2
\end{pmatrix}, \, F:= I_2 - E = \begin{pmatrix}
1/2 & -1/2\\ -1/2 & 1/2
\end{pmatrix}, \text{and then}$$ 
$$ 0 \leq Q = 
(P+A)\otimes E + (P-A)\otimes F,$$ so $P+A \geq 0 $
and $P-A \geq 0,$ and half of the proof is complete. The proof in the other direction is obtained by reversing the latest steps after supposing that we have a positive matrix 
 $P$ with $\|\mathrm{diag}(P)\|_\infty \leq 1 ,$ 
such that $P+A \geq 0 $ and $P-A \geq 0.$ 
\end{proof}

The theorem has the following immediate corollary.
\begin{corollary} \label{diff}
Let $A$ be a self-adjoint matrix in $M_n(\bc)$ then $$\|A\|_S = \min\{\|Q+R\|_S\,:\, Q,\, R\in M_n(\bc)_+,   \,  A  = Q-R\,\}.$$
\end{corollary}

\begin{proof}
If $A$ is a difference of two positive matrices, $A = Q-R,$ then for $P:= Q+R$ we have $-P \leq A \leq P, $ so by the theorem,  the Schur multiplier norm of $A$ is dominated by the norm of the diagonal of $P$ which is the Schur multiplier norm of $Q+R.$

By the theorem there exists a positive matrix $P$ such that $\|A|_S = \|\mathrm{diag}(P)\|_\infty = \|P\|_S$ and $-P \leq A \leq P.$ Then we can define positive matrices  $Q:= (P+A)/2$ and $R:= (P-A)/2,$ and they will satisfy the equations $P= Q+R$ and $A=Q-R,$ and the corollary follows.  
\end{proof}
\section{Mathias' construction of some norms and their dual norms on $M_{(m,n)}(\bc)$}
While working on the first version of this article we were not aware of Mathias' work \cite{Ma1}, which describes a way to determine the dual norms of certain norms on $M_{(m,n)}(\bc). $ One of the referees directed us to Mathias' work, and we do appreciate this help very much.
 We think Mathias' point of view on some pairs of dual norms clarifies several points in our research.

 Mathias' construction of dual norms  is clearly
inspired by the way Haagerup's result from Theorem \ref{Haa} is formulated. It is not difficult to see, that if we define a closed convex subset $\cas$ of the self-adjoint matrices in $M_{(m+n)}(\bc) $ by \begin{equation}
\cas:= \{ P_1 \oplus P_2\, : \, P_1 \in M_m(\bc)_+, \, P_2 \in M_n(\bc)_+ \text{ and } \|\mathrm{diag}(P_i)\|_\infty \leq 1 \}
\end{equation} then the function $N_{\cas}$ from  Theorem 4.1 of \cite{Ma1} becomes a norm, and Haagerup's result, Theorem \ref{Haa} of this article, tells that this norm will be the Schur multiplier norm. In item (c) of Corollary 4.3 in \cite{Ma1} Mathias states that the dual norm to the Schur multiplier norm on $M_{(m,n)}(\bc)$ may be obtained  from his Theorem 4.1, when this theorem is applied to the closed convex set $\cas^*$ defined as 

\begin{align*} \cas^* := \{ &D_1 \oplus D_2\, : D_1 \in M_m(\bc)_+, \, D_2 \in M_n(\bc)_+ ,\\& \, D_i \text{ diagonal and }  \mathrm{Tr}_m(D_1) \leq 1, \, \mathrm{Tr}_n(D_2) \leq 1\}.
\end{align*} 
It is quite easy to check that Mathias' Theorem 4.2 may be used to verify this statement. On the other hand we will like to investigate Mathias' norm $N_{\cas^*}$ based on his definition of this norm. The reason being that this will show us that for a complex $m \times n$ matrix $A$ we have   $N_{\cas^*}(A) = \|B_A\|_{cb}$ as defined in Theorem 2.4 of \cite{C3}. This norm is the completely bounded norm of the bilinear form with kernel $A$ on $(\bc^m,\|.\|_\infty)) \times  ( \bc^n, \|.\|_\infty),$ and this is the reason why we denoted this norm by $\|.\|_{cbB} $ in the works \cite{C3, C4}. In equation (3.5) from Theorem 3.6 of \cite{C3} it is stated that the norm $\|.\|_{cbB} $ is dual to the Schur multiplier norm. Here we will use  Mathias' Theorem 4.1 to recover the  part of the Theorem of 2.4 of \cite{C3} which describes the norm $\|.\|_{cbB}$ as a factorization norm
\begin{align}
\label{cbb}\forall A \in M_{(m,n)}(\bc):& N_{S^*}(A) = \|A\|_{cbB} \\ \notag  \|A\|_{cbB} = \min \{\|B\|_\infty :&  \,B \in M_{(m,n)}(\bc),\,  \exists \eta \in \bc^m, \, \eta_i \geq 0,\,\|\eta\|_2 = 1, \\ \notag \, \,  &\exists \xi \in \bc^n, \, \xi_j  \geq 0,\, \|\xi\|_2 = 1,\, A = \Delta(\eta)B\Delta(\xi)\}.
\end{align}
 Let $A$ be a complex $m \times n$ matrix, then  $N_{\cas^*}(A)  \leq 1$ if and only if there exist positive diagonal matrices  $D_1$ in $M_m(\bc)_+$ and $D_2$ in $M_n(\bc)_+$ both with trace at most $1$  such that \begin{equation} \label{NDA} \begin{pmatrix}
D_1& A\\A^* & D_2
\end{pmatrix} \geq 0.\end{equation}
Standard arguments will show that if some rows  or some  columns  in $A$  vanish then we may reduce the number number of columns or rows accordingly such that we only have to work with matrices $A$ where all columns and all rows are non trivial. The aspect of vanishing rows or columns  has previously  been treated in more details in both \cite{C3} and \cite{C4}.   It is easy to check that the value $N_{\cas^*}(A)$ is unchanged if we restrict the positive diagonal matrices $D_i$ to have the property Tr$_m(D_1) = 1,$ and Tr$_n(D_2) =1. $ 
If an  entry $d_{[1,(i,i)]} $ in $D_1$ vanishes, then the positivity of the block matrix from \eqref{NDA} implies that the entire row with number $i$ of $A$ vanishes. Since we do assume that $A$ has no trivial rows, all the $d_{[1,(i,i)]} > 0. $ A similar argument applies to the entries of $D_2$ and the columns of $A,$  so we may  assume that both $D_i$ are invertible.  
Define  positive unit vectors $\eta$ in $\bc^m$ and $\xi$ in $\bc^n$ by $\eta_i := \sqrt{d_{[1,(i,i)]}}$ and $\xi_j := \sqrt{d_{[2,(j,j)]} },$ so $D_1 = \Delta(\eta)^2$ and $D_2 = \Delta(\xi)^2, $ each $\eta_i >0$ and each $\xi_j > 0. $ Then we may elaborate on equation \eqref{NDA} as follows.
\begin{align} \label{NDA2} &0 \leq \begin{pmatrix}
\Delta(\eta)^2& A\\A^* & \Delta(\xi)^2
\end{pmatrix} = \\ &\notag   \begin{pmatrix}
\Delta(\eta)& 0\\ 0 & \Delta(\xi)
\end{pmatrix} \begin{pmatrix}
I_m & \Delta(\eta)^{-1} A\Delta(\xi)^{-1} \\   \Delta(\xi)^{-1} A^* \Delta(\eta)^{-1}  & I_n
\end{pmatrix} \begin{pmatrix}
\Delta(\eta_)& 0\\ 0 & \Delta(\xi)
\end{pmatrix}. 
\end{align} It is a well known fact that a block matrix of the  form $\begin{pmatrix}
I_m &B \\B^*&I_n
\end{pmatrix}$ is positive if and only if $\|B\|_\infty \leq 1, $ and it follows that the norm $N_{\cas^*}$ has the factorization property which characterizes the norm $\|.\|_{cbB}.$ We will extract the following proposition from the computations above.
\begin{proposition} \label{DualNormUnitVect} 
Given natural numbers $m$ and $ n,$ then the dual norm $\|.\|_{\cas^*} $ to the Schur multiplier norm $\|.\|_S$ on $M_{(m,n)}(\bc) $ is the norm $\|.\|_{cbB}.$ For any $A$ in $M_{(m,n)}(\bc)$ we have $\|A\|_{\cas^*} \leq 1$ if and only if there exist non negative unit vectors $\eta$ in $\bc^m$ and $\xi$ in $\bc^n$ such that \begin{equation} \label{PosMatrix}  \begin{pmatrix} \Delta(\eta)^2 & A \\ A^* &\Delta(\xi)^2 \end{pmatrix} \geq 0. \end{equation}
\end{proposition} 

The following proposition is based on the equality $\|.\|_{\cas^*} = \|.\|_{cbB}$ and \cite{C4}, Theorem 2.1 plus its Corollary 2.2. First we recall, that a complex $m \times n$ matrix in a natural way represents a linear map $L_A$ from the Hilbert space $(\bc^n, \|. \|_2) $ to the Hilbert space $(\bc^m, \|.\|_2).$ 

\begin{proposition} \label{unique} Let $A$ be a non zero  complex $m \times n $ matrix,  then  $A$ has a factorization $A = \Delta(\eta) B \Delta(\xi) $ such that $\eta_i \geq 0, $ $\|\eta\|_2 = 1,$ $\xi_j \geq 0,$ $ \|\xi\|_2 = 1,$ and $\|B\|_\infty = \|A\|_{\cas^*}. $ The positive unit vectors $\eta$ and $\xi$ are uniquely determined, and  $B$ is uniquely determined if the kernel of $L_B$ contains the kernel of $L_{\Delta(\xi)}$ and the kernel of $L_{B^*} $ contains the kernel of $L_{\Delta(\eta)}.$  

If $A$ is self-adjoint then $\eta = \xi$ and the uniquely determined factor  $B$ is self-adjoint. If $A$ is positive then  so is the uniquely determined  $B.$   
\end{proposition}
\begin{proof}
The uniqueness of the vectors $\xi$ and $\eta$ does 
does not follow directly from \cite{C4}. In that article the uniqueness is first established in the case where all the rows and all the columns in $A$ are non vanishing. The arguments right after equation  \eqref{NDA} tell that the Schur norm and also its dual norm for $A$ both equal the same norms for the {\em non vanishing  matrix} $\hat A, $ which is obtained from $A$ by deleting all trivial rows and columns.   We may then apply the results on unique factorization from \cite{C4} to  the matrix $\hat A$ defined on the set of indices, say $\cai$ and $\cj,$  with non vanishing rows and columns. In this way we get  we get a pair  of uniquely determined  {\em unit  vectors}  $\hat \xi$ and $\hat \eta$ with indices in  $\cai $ and $ \cj,$ such that the proposition holds for this reduced matrix.  We can then  define {\em unit vectors}  $\xi$ and $\eta$ in $\bc^n$ by supplementing the  vectors $\hat \xi, $ and $\hat \eta$ with a number of vanishing entries. Since the vectors $\hat \xi$ and $\hat \eta$ are uniquely determined as {\em unit vectors,} we can see that the {\em unit vectors } $\xi$ and $\eta$ are uniquely determined for the non vanishing matrix $A.$ 

In the case when A vanishes there are no candidates for the unique unit vectors in a factorization, exactly as in the case of a polar decomposition of the zero matrix. In that case there does not exist a  unique partial isometry part defined on the support-space of the matrix.    
\end{proof}

\section{The dual Schur norm of a self-adjoint matrix}

The computation of the norm $\|A\|_{\cas^*},$ for a self-adjoint matrix is based Proposition \ref{unique} and the methods used in the proof of  Theorem \ref{SchurNorm}.
\begin{theorem}
Let $A$ be a self-adjoint matrix in $M_n(\bc)$ then $\|A\|_{\cas^*} \leq  1$ if and only if there exists a vector $\g$ in $\br^n $ such that each $\g_j \geq 0,$  $\g_1 + \dots + \g_n =1 $ and $-\Delta(\g) \leq A \leq \Delta(\g).$ 
\end{theorem}  
 \begin{proof}
 Suppose $\|A\|_{\cas^*} \leq 1$ then by Proposition 
 \ref{DualNormUnitVect} there exist non negative unit vectors $\eta$ and $\xi$ in $\br^n$  such that \eqref{PosMatrix} is valid.  Define a non negative vector $\g$ in $\br^n$ by $\g_j := (\eta_j^2 +\xi_j^2)/2,$ then we have $\g_1 + \dots + \g_n =1$ and we can define a positive block matrix  $R$ by  
 $$R:=\begin{pmatrix} \Delta(\g)& A\\A& \Delta(\g) \end{pmatrix} = \frac{1}{2}\begin{pmatrix} \Delta(\eta)^2 & A\\A& \Delta(\xi)^2 \end{pmatrix}+ \frac{1}{2} \begin{pmatrix} \Delta(\xi)^2& A\\A& \Delta(\eta)^2 \end{pmatrix} \geq 0.$$
 If we use the notation from the proof of Theorem \ref{SchurNorm} we get $$ 0 \leq R = 
(\Delta(\g) +A)\otimes E + (\Delta(\g)-A)\otimes F,$$ so $\Delta(\g) +A \geq 0 $
and $\Delta(\g)-A \geq 0,$ and we see that the positivity conditions are necessary. To see the sufficiency, you can follow the steps above backwards. 
 \end{proof}
Later we will use a reformulation of the theorem above, and we will state this version as a corollary. 
 
 \begin{corollary} \label{cbBSa}
 
 Let $A$ be a self-adjoint matrix in $M_n(\bc)$ then $$\|A\|_{\cas^*} = \inf \{\mathrm{Tr}_n(\Delta(\lambda)) \, : \, \lambda \in \br^n, \text{ all } \lambda_j \geq 0, \, - \Delta(\lambda) \leq A \leq \Delta(\lambda)\}.$$

If $A$ is non zero then there exists a unique non negative vector $\lambda$ in $\br^n$  such that $$-\Delta(\lambda) \leq A \leq \Delta(\lambda) \text{ and  } \mathrm{Tr}_n(\Delta(\lambda)) = \|A\|_{\cas^*}.$$
 \end{corollary}
 
 \begin{proof} The first half of the corollary follows directly from the theorem. That the infimum is a minimum  follows from Mathias' Theorem 4.1, and hence  there exists a non zero vector $\lambda$ with non negative entries such that $\lambda_1 + \dots + \lambda_n = \|A\|_{S^*}$ and $-\Delta(\lambda ) \leq A \leq \Delta(\lambda).$ The uniqueness of $\lambda$ follows from Proposition \ref{unique}, after the following remarks. Since $\begin{pmatrix} \Delta(\lambda) & A \\A& \Delta(\lambda) \end{pmatrix} $ is positive, we get for any pair of vectors $\xi, \, \eta$ in $\bc^n$ and any  real $t > 0$ that
 \begin{align*} t\langle \Delta(\lambda)\xi, \xi \rangle + \frac{1}{t}\langle  \Delta(\lambda)\eta,\eta\rangle &\geq 2|\langle A \xi, \eta\rangle| \text{ by  minimising over } t 
  \\  \|\Delta(\lambda)^{(1/2)}\xi\| \|\Delta(\lambda)^{(1/2)} \eta\| & \geq 2|\langle A \xi, \eta\rangle |.\end{align*}
  Then, by standard operator theory, there exists a contraction matrix $C$ in $M_n(\bc)$ such that $A= \Delta(\lambda)^{(1/2)}C\Delta(\lambda)^{(1/2)},$ and we may  define a unit vector  $\xi$ with non negative entries by $\xi_j:= (\lambda_j/\|A\|_{S^*})^{(1/2)}$ for $1 \leq j \leq n,$  and a self-adjoint matrix $B$ with norm at most $\|A \|_{S^*}$ by $B:= \|A\|_{S^*} C.$ Since $\|A\|_{S^*} = \|A\|_{cbB},$ we get that  $A = \Delta(\xi)B\Delta(\xi) $ and this must be the unique  cbB-factorization of $A,$ so $\lambda$ is uniquely determined. 
 \end{proof}  

\section{ An example and Walter's inequality} \label{Exa}
In Theorem \ref{SchurNorm} we consider a complex self-adjoint $ n \times n$ matrix $A$ and positive matrices $P$  of the same type such that $P \geq A$ and $P\geq -A.$ There is a natural candidate for such a $P$ namely $P:= |A|:= \sqrt{A^2}$ and one may wonder if $|A| $ might be the optimal one, or at least, in some nice way, be  related to the one with minimal norm of the diagonal coming from Theorem \ref{SchurNorm}. Unfortunately this is not so, and we will provide an example to show this. 

In \cite{Wa} Martin E. Walter shows that for a complex $m \times n $ matrix $X$ the Schur multiplier norm is bounded from above by $\|X\|_S \leq \sqrt{  \|\mathrm{diag}(|X|)\|_\infty\|\mathrm{diag}(|X^*|)\|_\infty}.$
We will show how this inequality may be deduced from  Theorem \ref{SchurNorm} , and when this setup is  combined with the example we present, then we find that for any natural number $n$ there exists  a complex self-adjoint   $2n \times 2n $ matrix $Y$ which satisfies  $$  \|\mathrm{diag}(|Y|)\|_\infty = \sqrt{n}\|Y\|_S.$$

To see how Walter's inequality follows from Theorem \ref{SchurNorm} for a complex $m \times n $ matrix $X$ you may  look at the self-adjoint $(m+n) \times (m+n)$ matrix $Y$ given as  $Y := \begin{pmatrix} 0 &X\\X^* & 0\end{pmatrix}. $ It is well known, that for a $2 \times 2 $ complex  block matrix $Z$ in the form $Z= \begin{pmatrix}
0&A\\B&0
\end{pmatrix}$ we have
$\|Z\|_\infty = \max\{\|A\|_\infty , \, \|B\|_\infty\},$ and from here, it follows that $\|Z\|_S = \max\{\|A\|_S, \, \|B\|_S\}.$  Hence for $Y$ we get $\|Y\|_S= \|X\|_S = \|X^*\|_S.$  Then  for any real $t > 0 $ we may define  $P_t  = \begin{pmatrix} t |X^*| &0\\0 & (1/t)|X| \end{pmatrix}$ and we get that \begin{align*}P_t &= \begin{pmatrix}
\sqrt{t} &0 \\ 0 & \sqrt{1/t} \end{pmatrix}|Y|\begin{pmatrix}
\sqrt{t} &0 \\ 0 & \sqrt{1/t} \end{pmatrix} \\& \geq   \begin{pmatrix}
\sqrt{t} &0 \\ 0 & \sqrt{1/t} \end{pmatrix}Y\begin{pmatrix}
\sqrt{t} &0 \\ 0 & \sqrt{1/t} \end{pmatrix} = Y\\ & \text{ and }\\  
 P_t &= \begin{pmatrix}
\sqrt{t} &0 \\ 0 & \sqrt{1/t} \end{pmatrix}|Y|\begin{pmatrix}
\sqrt{t} &0 \\ 0 & \sqrt{1/t} \end{pmatrix} \\ & \geq   \begin{pmatrix}
\sqrt{t} &0 \\ 0 & \sqrt{1/t} \end{pmatrix}(-Y)\begin{pmatrix}
\sqrt{t} &0 \\ 0 & \sqrt{1/t} \end{pmatrix} = -Y. \end{align*} 
Then for $t = \sqrt{\|\mathrm{diag}(|X|)\|_\infty/\|\mathrm{diag}(|X^*|)\|_\infty}$ we may apply Theorem \ref{SchurNorm} to $P_t$ in order to obtain Walter's inequality. 
\begin{example} \label{Ex} Let $n$ be a natural number, $\{e_{(i,j)}\, :\, 1 \leq i,j \leq n \} $ the canonical set of matrix units for $M_n(\bc)$ and  the matrix $X:= \sum_{1 \leq i \leq n} e_{(i,1)}.$ Define the self-adjoint matrix $Y$ in $M_2(M_n(\bc))$  by $Y= \begin{pmatrix} 
0 & X \\X^*&0
\end{pmatrix} $ and the positive $2n \times 2n $ matrices $P_t$  are  defined as above from $|Y|.$  Then
\begin{itemize}
\item[(i)] $\|X\|_S = \|Y\|_S = 1$ 
\item[(ii)] $\|\mathrm{diag}(|X|)\|_\infty = \sqrt{n}$ and $\|\mathrm{diag}(|X^*|)\|_\infty = 1/\sqrt{n}$
\item[(iii)] $\|\mathrm{diag}(|Y|)\|_\infty = \sqrt{n}$ 
\item[(iv)] $\|\mathrm{diag}(P_{\sqrt{n}})\|_\infty = 1.$  
 \end{itemize}
\end{example} 

\begin{proof} 

 (i) We proved above that $\|X\|_S = \|Y\|_S.$ With respect to the norm $\|X\|_S,$ we remark that for any complex  $m \times n$ matrix $A,$ we have that the Schur product $X\circ A $ may be described  as $X\circ A = Ae_{(1,1)},$  so $\|X\|_S = 1.$ 

(ii) We have $X^*X = n e_{(1,1)} $ so $|X| = \sqrt{n} e_{(1,1)} $ and the first part of item (ii) follows. Let $E : = \sum_{(1 \leq i,j \leq n)}e_{(i,j)}, $ then $(1/n)E$ is a self-adjoint rank 1 projection and $XX^* = E = n\big((1/n)E\big),$ so $|X^*| = (1/\sqrt{n})E$ and the second half of item (ii) follows. 
\smallskip

(iii) We find that $|Y| = \begin{pmatrix}
|X^*| &0 \\ 0& |X|
\end{pmatrix},$ and item (iii) then follows from item (ii).
\smallskip

(iv) This follows from the computations just made. 

\end{proof} 

When we look at the self-adjoint $2n \times  2n$  matrix $Y$ in the above example we find that $\|\mathrm{diag}(|Y|)\|_\infty/ \|Y \|_S = \sqrt{n}, $ so the matrix $|Y| $ is not always a good candidate for the optimal matrix $P$ from Corollary  \ref{SchurNorm}.
\smallskip

With respect to Walter's inequality we get for this $Y$ that $Y= Y^*$ so Walter's\ inequality gives $$1 = \|Y\|_S \leq \sqrt{\|\mathrm{diag}(|Y|)\|_\infty \|\mathrm{diag}(|Y^*|)\|_\infty}= \sqrt{n}.$$ 

On the other hand Walter's inequality becomes an equality when applied to $X$ since we have 
 $$1 = \|X\|_S \leq \sqrt{\|\mathrm{diag}(|X|)\|_\infty \|\mathrm{diag}(|X^*|)\|_\infty}= 1.$$ 
 It seems that a closer analysis may shed some more light on the problem of deciding when Walter's inequality yields  a strong upper bound for the Schur multiplier norm. R. Mathias \cite{Ma2} has shown that if a complex $n \times n $ matrix $X$  satisfies $$ \mathrm{diag} (|X|) =  \mathrm{diag} (|X^*| ) = (1/n) \mathrm{Tr}_n(|X|) I_n $$ then $\|X\|_S = \sqrt{\| \mathrm{diag}(|X|)\|_\infty\|\mathrm{diag}(|X^*|)\|_\infty},$  and we have equality in Walter's inequality for such a matrix $X.$ This aspect of Walter's inequality is  discussed in more details in Section 4 of Davidson's and Donsig's article \cite{DD}. 

\section{Computation of the  $\cas^*$-norm as an infinite linear program} \label{LP}

The Corollary \ref{cbBSa} is easily reformulated into a {\em formal infinite linear program,} whose optimal value is the dual Schur norm  of the given  self-adjoint $n \times n $  matrix $A$. This linear
   minimizing  program is not entirely in  the standard form, since the number of linear constraints is infinite. On the other hand the Corollary \ref{cbBSa}
 shows that the program has an optimal solution.  We will use part of the words used in {\em linear programming theory,} so any vector representing a set  of variables is called a solution, the words {\em feasible solution} means a vector of variables, which satisfies all the linear constraints and an optimal solution is a feasible solution, which solves the infinite linear program.  We can not rely on the classical duality theorem for linear programs, so we will  show that the formal dual of this program has  an optimal solution with the same optimal value as the primal program, and we will get such a solution, which only involves a finite number of non vanishing variables. We will call any solution with only a finite number of non vanishing variables a {\em finite solution.} 
 The optimal finite solution to the dual program, which we produce,  is  found as a consequence of the duality between the 2 norms $\|.\|_{\cas^*}$ and $\|.\|_S$ and the characterization of the latter norm given in Corollary \ref{diff}.
   
 Let now $\cas$ denote the set of unit vectors in $\bc^n,$ then the infinite linear program $(\cp)$ is defined as: 
    \begin{align*} (\cp):  &\min \sum_{j = 1 }^n  \l_j\\  &\text{under the constraints }\\ & \forall j:  \l_j \geq 0\\ & \forall \xi \in \cas: \sum_{j=1}^n  \l_j|\xi_j|^2 \geq \langle A\xi, \xi\rangle\\ & \forall \eta \in \cas: 
   \sum_{j=1}^n  \l_j|\eta_j|^2 \geq -\langle A\eta, \eta\rangle.
   \end{align*}
  The formal dual program $(P^*)$ becomes:
  
  \begin{align}\notag (\cp^*):  &\max  \, \bigg(\sum_{\xi \in \cas} \mu_\xi \langle A\xi, \xi\rangle -  \sum_{\eta \in \cas} \nu_\eta\langle A\eta, \eta \rangle \bigg)\\ \notag 
  &\text{under the constraints }\\ \notag & \forall \xi \in \cas : \mu_\xi \geq 0 , \, \forall  \eta \in \cas: \nu_\eta \geq 0\\ \label{feasible} & \forall j \in \{ 1, \dots, n \}: \sum_{\xi \in \cas} \mu_\xi |\xi_j|^2 + \sum_{\eta \in \cas} \nu_\eta |\eta_j|^2 \leq 1.
   \end{align}
   
   \begin{remark} \label{convergence}  The expression 
   \begin{equation} \label{FesSum1}\bigg(\sum_{\xi \in \cas} \mu_\xi \langle A\xi, \xi\rangle -  \sum_{\eta \in \cas} \nu_\eta\langle A\eta, \eta \rangle \bigg)
   \end{equation} 
may not always  be well defined, but if an infinite solution \newline $\big((\mu_\xi)_{( \xi \in \cas)},\,  (\nu_\eta)_{( \eta \in \cas)}\big)$ is feasible, then we will prove that each of the sums in \eqref{FesSum1} converge absolutely. To see  that we look at a finite feasible solution $\big((\mu_\xi)_{(\xi \in S_0)} , \, (\nu_\eta)_{(\nu \in S_1)}\big).  $ Then we look at the sum in \eqref{FesSum2}, just below,  and use the feasibilty condition \begin{equation} \label{FesSum2} \sum_{\xi \in \cas_0}\mu_\xi + \sum_{\eta  \in \cas_1} \nu_\eta  = \sum_{j=1}^n \bigg(\sum_{\xi \in \cas_0}\mu_\xi|\xi_j|^2  + \sum_{\eta  \in \cas_1} \nu_\eta |\eta_j|^2\bigg) \leq n.\end{equation}
For each $j $ in $\{1, \dots, n\}$ the  sum of sums in the feasibility condition is always uniquely determined -  it may be infinite -   since all entries are non negative, and the sum of sums equals the supremum over all finite partial sum of sums. Then for  an infinite feasible solution $\big((\mu_\xi)_{( \xi \in \cas)},\,  (\nu_\eta)_{( \eta \in \cas)}\big)$ and two finite subsets $\cas_0 $ and $\cas_1$ of $\cas,$  the finite  solution $\big((\mu_\xi)_{( \xi \in \cas_0)},\,  (\nu_\eta)_{( \eta \in \cas_1)}\big)$ is also feasible and we may deduce that the following sum of sums is a sum of convergent sums of non negative reals which satisfy 
\begin{equation} \label{FesSum3} \sum_{\xi \in \cas}\mu_\xi + \sum_{\eta  \in \cas} \nu_\eta  \leq n.\end{equation}

  Since the numbers $|\langle A\xi, \xi\rangle |$ and $|\langle A\eta, \eta\rangle|$ are all bounded by $\|A\|_\infty,$ we get - for the feasible solution $\big((\mu_\xi)_{( \xi \in \cas)},\,  (\nu_\eta)_{( \eta \in \cas)}\big)$  - that the in the difference of sums in  \eqref{FesSum1} each of the sums converge absolutely since we have the inequality   \begin{equation} \label{AbsConv} \sum_{\xi  \in \cas}\mu_\xi\|A\|_\infty  + \sum_{\eta \in \cas}\nu_\eta \|A\|_\infty \leq n \|A\|_\infty.\end{equation}

For a feasible solution the inequality  \eqref{AbsConv} shows that the value  in \eqref{FesSum1} may be approximated  arbitrarily well by a difference  of sums over a finite number of indices, so we see that the value of an infinite feasible solution may be approximated arbitrarily well by the value of a finite feasible solution. On the other hand we will show, that the optimal value may be  obtained from a finite feasible solution, so we do not have to work with infinite optimal solutions in the arguments to come.   
   \end{remark}

   We will use the result stated in  the coming lemma a couple of times in the proof of Theorem \ref{LpSol}, and in order to present this lemma we need a definition first. 
   \begin{definition} Let $\xi $ be a unit vector in $\bc^n,$ then $ P_\xi$ denotes the matrix of the  orthogonal projection onto the space $\bc\xi.$
  \end{definition} 
   
   \begin{lemma} \label{lemFeas} 
    Let  $\big( (\mu_\xi)_{(\xi \in \cas)}, \, (\nu_\eta)_{( \eta \in \cas)} \big)$ be a finite solution to $(\cp^*)$  and define a positive matrix $P$ by \begin{equation} \label{ZfromSol} P := \sum_{\xi \in \cas} \mu_\xi P_\xi + \sum_{\eta \in \cas }\nu_\eta P_\eta. \end{equation}
   This finite solution is feasible if and only if $\|P\|_S \leq 1.$ 
   \end{lemma} 
  \begin{proof}
  For any vector $\g$ in $\cas$ and any $j$ in $\{1, \dots, n\}$ we have $$ |\g_j|^2 = | \langle \d_j,\g \rangle |^2 = \|P_\g \d_j\|^2 = \langle P_\g \d_j, \d_j \rangle = (P_\g)_{(j,j)}.$$
  Let us fix an index $j$ then \begin{equation}
 \sum_{\xi \in \cas} \mu_\xi |\xi_j|^2 + \sum_{\eta \in \cas }\nu_\eta |\eta_j|^2 =   \sum_{\xi \in \cas} \mu_\xi (P_\xi)_{(j,j)}  + \sum_{\eta \in \cas }\nu_\eta (P_\eta)_{(j,j)}  = P_{(j,j)},  
\end{equation} and the lemma follows from Schur's Theorem.   
\end{proof}    
   \begin{theorem} \label{LpSol} 
Let $A$ be a non zero self-adjoint matrix in $M_n(\bc).$    The optimal value for the programs $(\cp)$ and $ (\cp^*) $ is $\|A\|_{\cas^*},$ and an optimal solution to $(\cp)$ may be read out of the unique  cbB-factorization of $A.$ 
   
    Let  $Y$ be a self-adjoint matrix in $M_n(\bc)$ such that $\|Y\|_S = 1$ and  $\mathrm{Tr}_n(YA) = \|A\|_{\cas^*}.$ 
    Then a finite  optimal solution for $(\cp^*)$ may be obtained from any decomposition of $Y,$ as described in Corollary \ref{diff}, $ Y  = Q-R$ with $Q$ and $R$ positive such that $\|Y\|_S = \|Q+R\|_S.$

Let $\big( (\mu_\xi)_{(\xi \in \cas)}, \, (\nu_\eta)_{( \eta \in \cas)} \big)$ be a finite optimal  solution to $(\cp^*),$ 
and let $Y$ be defined by 
\begin{equation} \label{DefQ}  Y := \sum_{\xi \in \cas} \mu_\xi P_\xi - \sum_{\eta \in \cas }\nu_\eta P_\eta. \end{equation} 
then $\|Y\|_S = 1 $ and $\mathrm{Tr}_n(YA) = \|A\|_{\cas^*}.$ 
   \end{theorem}
   \begin{proof}
   It is elementary to see that the program $(\cp)$ consists in finding the infimum $$\inf\{\mathrm{Tr}_n(\Delta(\l))\, : \, -\Delta(\l) \leq A \leq \Delta(\l)\}.$$
   According to Corollary \ref{cbBSa} this infimum is a minimal value equal to $\|A\|_{\cas^*}, $ which is obtained from the decomposition\newline  $A = \Delta(\xi)B\Delta(\xi ) $ such that an  optimal $\l$ is  given by  $\l_j := \|A\|_{\cas^*}\xi_j^2.$  
   \smallskip

Let us turn to the dual program, and let $Y$ be a self-adjoint matrix such that $\|Y\|_S =1$ and $\mathrm{Tr}_n(YA)= \|A\|_{\cas^*}.$ We will construct a finite  optimal solution for $(\cp^*)$ based on this matrix.
Before we start to compute, we remind the reader, that for two self-adjoint $n \times n$ matrices $A$ and $Y$ the  the value $\mathrm{Tr}_n(AY)$ is real since $\mathrm{Tr}_n(AY) = \mathrm{Tr}_n(YA) = \mathrm{Tr}_n((AY)^*) = \overline{\mathrm{Tr}_n(AY)}.$ Then  we are dealing with real numbers  in the computations to come, and we may use inequality signs  without  further comments. By Corollary \ref{diff} we can write $Y = Q - R$ with $Q$  and $R$ positive such that $\|Q+R\|_S = 1.$  We will first suppose that both of the matrices $Q$ and  $R$ do not vanish. Then  there exists an orthonormal set of eigenvectors $\{\xi_1, \dots \xi_k\},$   and strictly positive real numbers $ \{\mu_1 , \dots, \mu_k\}$ such that $Q=\sum_s \mu_s P_{\xi_s}.  $ In the same way  we can find an orthonormal  set of eigenvectors $\{\eta_1, \dots, \eta_l\},$  and strictly positive real numbers  $\{\nu_1, \dots, \nu_l\}$ such that  $R =\sum_t \nu_t P_{\eta_t} .$ 
    
   We will  define the coefficients $\mu_\xi$ and $\nu_\eta$ of a finite  solution to  $(\cp^*)$ by $$ \mu_\xi := \begin{cases} \mu_s \text{ if } \xi =\xi_s \\ 0 \,\, \text{ else} \end{cases} \quad  \nu_\eta := \begin{cases} \nu_t \text{ if } \eta  =\eta_t \\ 0 \,\, \text{ else} \end{cases}.$$ Let us show that this constitutes a finite feasible solution, and to do that we go back to equation \eqref{ZfromSol} and define $P$ from this equation and   since  $P= Q+R,$ we get by Schur's theorem that $\|\mathrm{diag}(P)\|_\infty = 1$ and by Lemma  \ref{lemFeas} we find that this solution is a finite feasible solution.  Next we  will  compute the value of this finite feasible solution.
 \begin{align}
 &  \sum_{\xi \in \cas} \mu_\xi \langle A\xi, \xi\rangle -  \sum_{\eta \in \cas} \nu_\eta\langle A\eta, \eta \rangle = \sum_{s=1}^k \mu_s\langle A\xi_s, \xi_s\rangle - \sum_{t = 1}^l \nu_t \langle A\eta_t, \eta_t \rangle \\ \notag =& \sum_{s=1}^k \mu_s\mathrm{Tr}_n( AP_{\xi_s} )  - \sum_{t = 1}^l \nu_t \mathrm{Tr}_n( A P_{\eta_t}) =    \mathrm{Tr}_n(A(Q-R ) )\\ \notag =&  \mathrm{Tr}_n(AY)  =   \|A\|_{\cas^*}.
\end{align}   

If the decomposition $A = Q-R$ comes with $R= 0$ in $M_n(\bc),$ then $A$   is positive. In this case we may proceed as above but define the solution to be given as 
\begin{equation} \label{R=0}  \mu_\xi := \begin{cases} \mu_s \text{ if } \xi =\xi_s \\ 0 \,\, \text{ else} \end{cases} \quad  \forall \eta \in \cas:\,  \nu_\eta := 0.\end{equation} 

If $A = Q- R $ comes with $Q= 0$ we define in the same way a solution based on $R = -A$ by 
\begin{equation}\label{Q=0}  \forall \xi \in \cas:\, \mu_\xi:= 0  \quad  \nu_\eta := \begin{cases} \nu_t \text{ if } \eta  =\eta_t \\ 0 \,\, \text{ else} \end{cases}.\end{equation} 

In both of the cases, $R =0$ and $Q \neq 0$ or $R \neq 0$ and $Q = 0,$ we may, in the same way as in the general case,  establish that the solutions from both \eqref{R=0}  and  \eqref{Q=0} are feasible and both have the value $\|A\|_{S^*}.$

These computatiions show that the optimal value for $(\cp^*)$ is at least $\|A\|_{\cas^*},$ and that the value $\|A\|_{\cas^*}$ may be obtained from a finite feasible solution. 

Since we do not know that  {\em the weak  duality result for linear programs}  holds in the case we consider, we will show that the optimal and maximal value for $(\cp^*)$ is at most $\|A\|_{\cas^*}.$ \newline 
Let  $\big((\xi_1, \mu_1), \dots (\xi_k, \mu_k)\big),(\eta_1, \nu_1), \dots, (\eta_l,\nu_l)\big)$  be a finite  feasible solution, then we will  show, that its value is at most  $\|A\|_{\cas^*}.$  We will define  positive matrices $Q$ and $R$ by $$ Q:= \sum_{s=1}^k \mu_s P_{\xi_s} , \quad R:= \sum_{t=1}^l \nu_t P_{\eta_t} ,$$ then by Lemma \ref{lemFeas}, the feasibility of the solution shows that $\|Q+R\|_S \leq 1.$  
We can then define a self-adjoint matrix $Y: = Q - R$ and and a positive matrix $P:= Q+R$ which satisfy $$ - P = -Q-R \leq Q-R = Y \leq Q+R = P,$$ so by Theorem \ref{SchurNorm}  we have $\|Y\|_S \leq \|P\|_S \leq 1.$
Let us now bound  the value of our feasible solution 
   \begin{align}
 &  \sum_{\xi \in \cas} \mu_\xi \langle A\xi, \xi\rangle -  \sum_{\eta \in \cas} \nu_\eta\langle A\eta, \eta \rangle = \sum_{s=1}^k \mu_s\langle A\xi_s, \xi_s\rangle - \sum_{t = 1}^l \nu_t \langle A \eta_t, \eta_t \rangle \\ \notag =& \sum_{s=1}^k \mu_s\mathrm{Tr}_n( AP_{\xi_s} )  - \sum_{t = 1}^l \nu_t \mathrm{Tr}_n( A P_{\eta_t}) =    \mathrm{Tr}_n(A(Q-R) ) = \mathrm{Tr}_n(AY) \\ 
 \notag \leq & \|Y\|_S\|A\|_{\cas^*} \leq \|A\|_{\cas^*}.
\end{align}  
By the Remark \ref{convergence} it follows that the value of any feasible solution may be approximated arbitrarily well by the value of a finite feasible solution so the weak duality holds, and the maximal value of $(\cp^*)$ must be $\|A\|_{\cas^*}.  $ We do also see, that a finite optimal solution for $(\cp^*)$ may be obtained from any self-adjoint matrix $Y$ with $\|Y\|_S = 1, $ which satisfies Tr$_n(YA) =\|A\|_{\cas^*}, $ so the second statement in the theorem follows.  

The third statement in the theorem follows from the arguments given in the proof of the second statement as follows. For a given finite optimal solution $\big((\xi_1, \mu_1), \dots, (\xi_k, \mu_k), ( \eta_1, \nu_1) , \dots, (\eta_l, \nu_l)\big)$ we define as above $$ Q := \sum_{s =1}^k \mu_sP_{\xi_s}, \quad R: = 
\sum_{t =1}^l\nu_tP_{\eta_t},\quad Y:= Q-R, \quad P:= Q+R.$$ Then the feasibility shows that for the positive matrix $P$ we have $\|P\|_S \leq 1$ and the inequalities $- P \leq Y \leq P $ show that $\|Y\|_S \leq \|P\|_S\leq 1. $ Finally the optimality shows that $\mathrm{Tr}_n(AY) = \|A\|_{\cas^*}$ so $\|Y\|_S = 1, $ and the last statement in the theorem follows. 
   \end{proof}

\end{document}